\documentclass[12pt]{amsart}
\usepackage{amsmath,amssymb,amscd}
\usepackage{epsfig}

\newtheorem{thm}{Theorem}[section]
\newtheorem{prop}[thm]{Proposition}
\newtheorem{coro}[thm]{Corollary}
\newtheorem{lemma}[thm]{Lemma}

\newtheorem{conj}[thm]{Conjecture}

\theoremstyle{definition}

\newtheorem{rmk}[thm]{Remark}
\newtheorem{rem}[thm]{Remark}
\newtheorem{eg}[thm]{Example}
\newtheorem{ex}[thm]{Exercise}
\newtheorem{defn}[thm]{Definition}
\newtheorem{notn}[thm]{Notation}

\newcommand\wt[1]{\widetilde{#1}}

\def\til{\widetilde}
\def\mbb{\mathbb}
\def\mcl{\mathcal}
\def\mfk{\mathfrak}

\def\ten{\otimes}
\def\ex{\times}
\def\strl{\stackrel}
\def\tu{\textup}

\def\lan{\langle}
\def\ran{\rangle}
\def\a{\alpha}
\def\b{\beta}
\def\d{\delta}
\def\D{\Delta}

\def\G{\Gamma}
\def\D{\Delta}

\def\l{\lambda}

\def\sm{\sigma}

\def\o{\omega}

\def\P{\mbb P}

\def\Q{\mbb Q}
\def\cT{\mathcal T}
\def\cO{\mathcal O}

\def\inj{\hookrightarrow}

\def\dra{\dashrightarrow}
\def\SS{\mcl{O}}

\def\rk{\textup{rank} \, }
\def\lra{\longrightarrow}

\def\spec{\tu{Spec }}
\def\proj{\tu{Proj }}

\def\Proj{\tu{\textbf{Proj}\,}}
\def\Spec{\tu{\textbf{Spec}\,}}
\def\mod{/ \! \! /}

\def\SL{\textup{SL}}
\def\om2{\omega^{\ten 2}}
\def\Gr{\tu{Grass}}
\def\Mg{\overline{M}_g}
\def\Mgn{\overline{M}_{g,n}}
\def\KMg{K_{\Mg}}

\def\Ch{{\bf Chow}}

\def\FMps{\overline{{\bf \mcl M}}^{ps}_g}
\def\Mps{\overline{M}^{ps}_g}
\def\dps{\delta^{ps}}
\def\FMg{\overline{\bf \mcl M}_g}
\def\FMgone{\overline{\bf \mcl M}_{g,1}}
\def\KFMg{K_{\FMg}}

\def\ds{\oplus}

\def\inj{\hookrightarrow}

\def\PGL{\tu{PGL}}

\def\Sym{\tu{Sym}}

\def\Gr{\mathrm{Grass}}
\def\FGr{\underline{\Gr}}
\def\ra{\rightarrow}
\def\pj{\sum_{j=1}^rp_j}
\def\t{\tau}

\def\bar{\overline}

\input xy
\xyoption{all}
\input epsf
\newcommand{\bP}{{\mathbb P}}

\newcommand{\bG}{{\mathbb G}}

\newcommand{\bQ}{{\mathbb Q}}

\newcommand{\cB}{{\mathcal B}}

\newcommand{\cD}{{\mathcal D}}
\newcommand{\cE}{{\mathcal E}}
\newcommand{\cF}{{\mathcal F}}

\newcommand{\cU}{{\mathcal U}}
\newcommand{\ocM}{\overline{\mathcal{M}}}
\newcommand{\cX}{{\mathcal X}}
\newcommand{\cY}{{\mathcal Y}}
\newcommand{\cZ}{{\mathcal Z}}
\newcommand{\Ex}{\mathrm{Ex}}

\title[Log canonical models]{Log canonical models
for the moduli space of curves: First divisorial
contraction}
\author{Brendan Hassett}
\address{Rice University, 6100 Main St., Houston TX 77251-1892}
\author{Donghoon Hyeon}
\address{Northern Illinois University,  DeKalb IL 60115
}
\date{\today}
\begin{document}
\begin{abstract}
In this paper, we initiate our investigation of log canonical models
for $(\bar{\mcl M}_g,\a \delta)$ as we decrease $\a$ from 1 to 0. We
prove that for the first critical value $\a = 9/11$,  the log
canonical model is isomorphic to the moduli space of pseudostable
curves, which have nodes and cusps as singularities. We also show
that $\a = 7/10$ is the next critical value, i.e., the log canonical
model stays the same in the interval $(7/10, 9/11]$. In the
appendix, we develop a theory of log canonical models of stacks that
explains how these can be expressed in terms of the coarse moduli
space.
\end{abstract}
\maketitle

\tableofcontents
\section{Introduction}
This article is the first of a series of papers investigating
the birational geometry of the moduli space
$\Mg$ of stable curves of genus $g$.  The guiding problem
is to understand the canonical model of $\Mg$
\[
\proj \left(\ds_{n\geq 0} \Gamma(\Mg, n K_{\Mg})\right),
\]
when $\Mg$ is of general type.
The moduli space is known to be of general type for
$g\ge 24$ by work of Eisenbud, Harris,
and Mumford \cite{HM};  recently, Farkas has announced results
for $g=22$ and $23$.

One significant obstacle is that
the canonical ring of $\Mg$ is not known to be
finitely generated for {\em any} genus $g$.
A fundamental conjecture of birational geometry predicts that
canonical rings of varieties of general type are finitely generated,
but this has been proven only in small dimensions.
To get around this difficulty, we
study the intermediate {\it log canonical models}
\[
\Mg(\a):= \proj \left(\ds_{n\geq 0} \Gamma(n (\KFMg +
\a\d))\right)
\]
for various values of $\a \in \Q \cap [0,1]$.  Here $\FMg$ is the moduli stack,
$\d_i,i=0,\ldots, \lfloor g/2 \rfloor$ its boundary divisors,
and $\KFMg$ its canonical divisor.  Let $\FMg \ra \Mg$ denote the morphism
to the coarse moduli space, which has
boundary divisors $\D_i$ and canonical class $\KMg$.

A general discussion of log canonical models for stacks is offered
in Appendix~\ref{sect:LCMS}, but we sketch the basic idea here.
Assume that $g\ge 4$;  the $g=2$ case is discussed in \cite{Has}.
Then the $\bQ$-divisor on the coarse moduli space
\[
K_{\Mg}+\alpha(\D_0+\D_2+\ldots+\D_{\lfloor g/2\rfloor})+
\frac{1+\alpha}{2}\D_1  \tag{$\dag$}
\]
pulls back to $\KFMg+\a \d$ on the moduli stack
(see Remark 4.9 of \cite{Has}).  We shall often
abuse notation by conflating these divisors.
The universal property of the coarse moduli space implies that sections
of invertible sheaves on $\FMg$ are all pull-backs of
sections of the corresponding reflexive sheaves on
$\Mg$.  The log canonical model of $\FMg$ with respect to
$\KFMg+\a \d$ can thus be identified with the log canonical model
of $\Mg$ with respect to $(\dagger)$ (see Proposition~\ref{prop:redcoarse}).

For large values of $\a$ these log canonical models are well-understood.
Results of Mumford, Cornalba, and Harris imply that $\Mg(\a)\simeq \Mg$
if and only if $9/11<\a\le 1$.  Indeed, by the uniqueness
of log canonical models, it suffices to show the
pair $(\FMg,\a\d)$ has log canonical singularities and
$\KFMg+\a\d$ is ample on $\Mg$,
for $\a$ in this range.   By Theorem 2 of \cite{HM},
$\Mg\setminus \Delta_1$ actually has canonical
singularities;  the pair $(\FMg,\a\d)$ is a smooth stack
with a normal-crossings divisor and thus is log canonical.
The main result of \cite{CH} says that
$$a\lambda-\delta=\frac{a}{13}(\KFMg+(2-\frac{13}{a})\delta)$$
is ample for $a>11$.  Indeed, $(\FMg,\a\d)$ is a
{\em strict log canonical model } in the sense that all
divisors with negative discrepancies lie over $\d$
(cf. \cite{Has} \S 4.3).

In this paper, we describe what happens when $\a = 9/11$.
The pair $(\FMg,9/11\ \d)$ remains log terminal, but
$\KFMg+9/11\ \d$ is not ample on $\Mg$.  Mumford \cite{M} described
an extremal ray $R$ for this divisor, parametrizing the
\textit{elliptic tails}.  Recall that the normalization of
$\delta_1$ is isomorphic to $ \overline{\bf \mcl M}_{1,1}\times
\overline{\bf \mcl M}_{g-1,1} $ and $\overline{M}_{1,1}\simeq
\bP^1$. Our ray is the class of the fiber under projection to the
second factor.
\begin{figure}[ht]\label{F:et}
\centerline{\psfig{figure=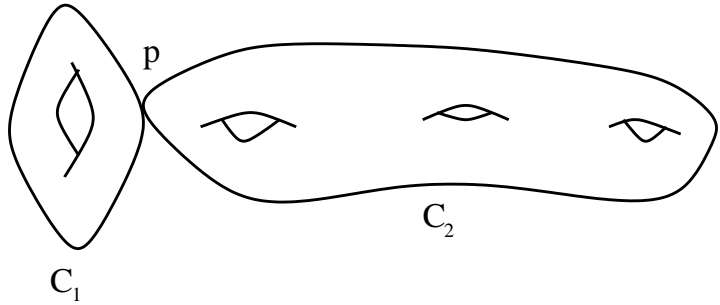}} \caption{Fix $C_2$ and
vary the elliptic curve}
\end{figure}
Explicitly, for each $(C_2,p)$ of $\overline{M}_{g-1,1}$
consider the stable curves
$C=C_1 \cup_p C_2,$
where $(C_1,p)$ is a varying curve in $\overline{M}_{1,1}\simeq \P^1$,
meeting $C_2$ at a single node $p$.  These are parametrized by
a rational curve $R(C_2,p)$;  the class $R=[R(C_2,p)]$
is independent of $C_2$ and $p$ and
$$(\KFMg+9/11\ \d).R=0.$$
This is the {\em unique} extremal ray for $\KFMg+9/11\ \d$;
it meets every other curve class in $\Mg$
positively (by the analysis of \cite{CH} or \S 6 of \cite{GKM}).

Applying the basepoint-freeness theorem to $\KFMg+9/11\ \d$
(see Corollary \ref{coro:bpfs}), we obtain an extremal contraction
of $R$
$$\Upsilon: \Mg \ra \Mg(9/11).$$
Since $R\cdot \D_1<0$, $\Upsilon$ is birational with exceptional locus
contained in $\D_1$.  More precisely, its restriction
to the generic point of $\D_1$ is
projection onto the $\overline{M}_{g-1,1}$ factor.
Thus $\Upsilon$ is a divisorial contraction, and $\Mg(9/11)$ remains
$\Q$-factorial.

Is $\Mg(9/11)$ the coarse
moduli space of some moduli stack of curves with prescribed properties?
Such a moduli functor is implicit in the work of D. Schubert
\cite{Schubert}. He constructed the moduli space
$\Mps$  of pseudostable curves applying geometric invariant
theory (GIT) to tricanonical curves.  Recall that a {\em cusp}
is a curve singularity analytically-locally isomorphic
to $y^2=x^3$.  A complete curve is
{\em pseudostable} if
\begin{enumerate}
\item it is connected, reduced, and
has only nodes and cusps as singularities;
 \item every subcurve of genus one meets the rest of
the curve in at least two points;
\item the canonical sheaf of the curve is ample.
\end{enumerate}
The last condition means that each subcurve of genus zero
meets the rest of the curve in at least three points.
The corresponding moduli stack is denoted $\FMps$.
\begin{thm}\label{T:main1}
There is a morphism of stacks
\[
\cT : \FMg \to \FMps
\]
which is an isomorphism in the complement of $\delta_1$.
For a stable curve $C\in \delta_1$, $\cT(C)$ is obtained
by replacing each elliptic tail of $C$ with a cusp.
$\Mps$ and $\Mg(9/11)$ are isomorphic compactifications
of the moduli space of curves and
the induced morphism on coarse moduli spaces
\[
T: \Mg \ra \Mps
\]
coincides with the extremal contraction $\Upsilon$.
\end{thm}
See \S \ref{S:nattran} for a
precise formulation of what it means to replace an elliptic
tail with a cusp.
We can characterize the $\a$ for which $\Mg(\a)$ is
isomorphic to the moduli space of pseudostable curves:
\begin{thm}\label{T:main2}
For $7/10<\a\le 9/11$, $\Mg(\a)$ exists as a projective variety
and is isomorphic to $\Mps$.
\end{thm}
What happens to $\Mg(\a)$ at the next
critical value $\a = 7/10$?  Here we obtain a \textit{small
contraction} (see Remark~\ref{rmk:codimtwo})
and one naturally seeks a construction of its
logarithmic flip.  In a subsequent paper,
we shall show that $\Mg(7/10)$ is isomorphic to the GIT quotient
of the Chow variety of bicanonical curves and  its flip is
a GIT quotient of the Hilbert scheme of bicanonical curves.

Finally, we survey recent results of a similar flavor.
Shepherd-Barron has analyzed the canonical model of the moduli
space of principally polarized abelian varieties of dimension
$\ge 12$ \cite{SB} Theorem 0.2.  Keel and Tevelev have
considered log canonical models for moduli of hyperplane
arrangements \cite{KT} 1.14.  These can be understood explicitly
in special cases, e.g., $(\overline{M}_{0,n},\d)$ is a log canonical
model \cite{HasW} 7.1.

Throughout, we work over an algebraically closed field $k$ of
characteristic zero.

\subsubsection*{Acknowledgments} The first author was partially
supported by National Science Foundation grants 0196187 and 0134259,
the Sloan Foundation,     and the Institute of Mathematical Sciences
of the Chinese University of Hong Kong. The second author was
partially supported by the Korea Institute for Advanced Study.  We
owe a great deal to S. Keel, who helped shape our understanding of
the birational geometry of $\Mg$ through detailed correspondence. We
are also grateful to Y. Kawamata and B.P. Purnaprajna for useful
conversations.

\section{Resum\'e of Schubert's thesis}\label{S:schubert}
In this section, we sketch Schubert's construction of the moduli
space of pseudostable curves \cite{Schubert}.  Throughout this
section we assume $g > 2$.  (The special features of the $g=2$ case
are addressed in \cite{HL}.)

\begin{defn}
An {\em $n$-canonical curve}
$$C\subset \P^N, \quad  N:=\begin{cases} g & \text{ if }n=1 \\
                (2n-1)(g-1)-1 & \text{ if } n>1
                \end{cases}
$$
is a Deligne-Mumford stable curve imbedded by the sections
of the pluricanonical bundle $\omega_C^{\otimes n}$.
($\omega_C^{\otimes n}$ is very ample
for $n\ge 3$ \cite{DM}, Theorem 1.2.)
\end{defn}
These form a locally closed
subset of the Chow variety (or Hilbert scheme) parametrizing
curves in $\P^{(2n-1)(g-1)-1}$ of degree $d:=2n(g-1)$.

Mumford gave the following GIT construction of $\Mg$ \cite{M}:
Fix $n\ge 5$ and consider the
closure $\Ch_n$ of $n$-canonical curves in the Chow variety,
with the natural induced action of $\SL_{N+1}$.
The stable points (in the sense of GIT) are precisely the
$n$-canonical curves, hence
$$\Ch_n \mod \SL_{(2n-1)(g-1)} \simeq \Mg, \quad n\ge 5.$$

Schubert investigated the stability of points in $\Ch_3$.
He proves that a Chow point is GIT stable if and only
if the corresponding cycle is a reduced pseudostable curve
of genus $g$.  Furthermore, there are no strictly semistable
points.  Here is the idea of the argument:
By general results from \cite{M}, smooth
tricanonical curves have GIT-stable Chow points.  Using
the methods of \cite{M} for finding destabilizing one-parameter
subgroups, one checks that the cycles associated
to the following types of curves are necessarily
unstable:
\begin{enumerate}
\item{curves with singularities other than nodes or cusps,
including nonreduced curves;}
\item{curves $C\subset \P^{5g-6}$ with
$\SS_C(+1)\not \simeq \omega_C^{\otimes 3}$;}
\item{curves contained in a proper linear subspace
of $\P^{5g-6}$;}
\item{curves containing an elliptic tail, i.e., curves
containing a subcurve of genus one meeting the rest of
the curve in a single point.}
\end{enumerate}
By elimination, only tricanonically-embedded
pseudostable curves have a chance
of having stable Chow points.

The proof that these do have stable Chow points
is indirect:  Let $C\subset \P^{5g-6}$ be a tricanonically-embedded
pseudostable curve.
A deformation-theoretic argument
shows there a exists a flat family of curves in $\P^{5g-6}$
$$\pi: \mcl C \ra \D=\spec k[[t]]$$
with closed fiber $C$ and generic fiber smooth and tricanonically
embedded.  Consider the induced morphism to the GIT quotient
$$\mu: \D \ra \Ch_3\mod \SL_{5g-5};$$
of course, $\mu(0)$ must correspond to certain semistable points in $\Ch_3$.
After base change $\widetilde{\D}\ra \D$ (substitute $t=\tau^M$ for some $M$),
we obtain a lift
$$\widetilde{\mu}:\widetilde{\D} \ra \Ch_3^{ss}$$
into the semistable locus of the Chow variety.  The generic
point maps into the stable locus.
Let ${\mcl D}\ra \widetilde{\D}$ be the family of pseudostable
curves associated to $\widetilde{\mu}$.  Note that $\mcl D$
and the pull-back $\widetilde{\mcl C}:=\mcl C\times_{\D}\widetilde{\D}$
agree up to the action of $\SL_{5(g-1)}$.
To prove that $C$ is semistable, we
show that the closed fibers of $\mcl C$ and $\mcl D$ agree.
We require the following uniqueness result, which may be regarded the valuative
criterion for separatedness for the moduli functor of pseudostable
curves:
\begin{prop}[cf. Lemma 4.2 of \cite{Schubert}]
Let ${\mcl D}_1,{\mcl D}_2\ra \D$ be flat families of pseudostable
curves with smooth generic fibers,  isomorphic over $k((t))$.
Then the closed the fibers are isomorphic as well.
\end{prop}
Finally, suppose there there were strictly semistable
points:  Then some point of $\Ch_3^{ss}$ would have
positive-dimensional stabilizer, but pseudostable curves have
finite automorphism groups \cite{Schubert}, pp 312.

The following moduli functor is implicit is Schubert's work:
\begin{defn} The moduli stack of pseudostable curves $\FMps$
is defined
\[
\FMps(S) =  \left\{ f : C \to S \, \, \left|  \begin{array}{l}   \,
\mbox{
(i) \, $f$ is a proper and flat morphism.} \\
 \, \mbox{ (ii) \, geometric fibres of $f$ are pseudo-} \\
 \, \mbox{\quad stable curves of genus $g$.}
\end{array}\right.
\right\}
\]
where $S$ is a scheme of finite type over $k$.
\end{defn}
This is a Deligne-Mumford stack, as pseudostable curves
have finite automorphism groups when $g>2$.
Its coarse moduli space is denoted $\Mps$.
\begin{thm}\label{T:Mps}(cf. \cite{Schubert}, Theorem 5.4)
The quotient stack $[\Ch_3^s/\PGL_{5g-5}]$ is isomorphic
to $\FMps$ and
$$\Ch_3 \mod \SL(5g-5) \simeq \Mps.$$
\end{thm}
\begin{proof}
Since $\Ch^s_3$ parametrizes
a family of pseudostable curves,
we have a morphism of stacks
$$\iota:[\Ch^s_3/\PGL_{5g-5}] \ra \FMps.$$
Two tricanonically-embedded pseudostable curves are
isomorphic if and only if they are projectively equivalent,
so $\iota$ is an open embedding.  The quasiprojective
variety $\Ch^s_3 \mod \SL_{5g-5}$ is the coarse moduli space
for $[\Ch^s_3/\PGL_{5g-5}]$ by \cite{GIT} Theorem 1.10.
But $\Ch_3$ has
no strictly semistable points under the $\SL_{5g-5}$ action, so
$\Ch^s_3 \mod \SL_{5g-5}=\Ch_3 \mod \SL_{5g-5}$.
Since $[\Ch^s_3/ \PGL_{5g-5}]$ has a projective coarse moduli
space, it is necessarily proper, and
$\iota$ is an isomorphism.
\end{proof}

\begin{rmk}
Let $C$ be a stable curve in $\delta_1$:  Schubert's reasoning
proves the Chow point of
the {\em four-canonical} image of $C$ is also GIT unstable.
Indeed, all the stability analysis sketched above applies to
four-canonically embedded pseudostable curves.
Consequently, we see that
$$\Mps\simeq \Ch_4\mod \SL_{7g-7}.$$
\end{rmk}

\section{Constructing the natural transformation $\cT$}
\label{S:nattran}
Here we prove the first part of
Theorem~\ref{T:main1}.  Throughout, we assume $g>2$.

\subsection{Replacing elliptic tails with cusps}
\begin{prop}   \label{prop:replacing}
Let $C$ be a stable
curve of genus $g>2$ with elliptic tails $E_1,\ldots,E_r$, i.e.,
connected genus-one subcurves meeting
the rest of $C$ in a single
node.  Let $D$ denote the union of the components of $C$ other than
elliptic tails and $p_i$ the node where $E_i$ meets $D$,
$i=1,\ldots,r$. Then there exists a unique curve $\cT(C)$
characterized by the following properties:
\begin{enumerate}
\item{there is a birational morphism $\nu:D\ra \cT(C)$,
which is an isomorphism away from $p_1,\ldots,p_r$;}
\item{$\nu$ is bijective and maps each $p_i\in D$
to a cusp $q_i\in \cT(C)$, formally isomorphic
to $y^2=x^3$.}
\end{enumerate}
There is a unique replacement morphism $\xi_C:C\ra \cT(C)$ with
$\xi_C|_D=\nu$ and $\xi_C|_{E_i}$ is constant. Note that $\cT(C)$
has arithmetic genus $g$.
\end{prop}
\begin{proof}
To determine $\cT(C)$, it suffices to specify the subrings
$$\SS_{\cT(C),q_i}\subset \SS_{D,p_i}.$$
Let ${\mfk m}_{D,p_i} \subset \SS_{D,p_i}$ be the maximal ideal and define
$\SS_{\cT(C),q_i}$ to be the algebra generated by the constants
and ${\mfk m}^2_{D,p_i}$.
A local computation shows that $\cT(C)$ has a cusp
at $q_i$.  Indeed, ${\mfk m}_{\cT(C),q_i}$
is generated by two elements $\hat{x}$ and $\hat{y}$,
vanishing at $p_i$ to orders two
and three respectively, whence $\hat{y}^2=c\hat{x}^3$ for
some unit $c\in \SS_{D,p_i}$.    Conversely, any germ of a cuspidal curve
normalized by $(D,p_i)$
is obtained in this way.  The morphism $\nu$ is the normalization
of the cusps $q_1,\ldots,q_r$.
\end{proof}

How do we describe a morphism of stacks
$$\cT:\FMg \ra \FMps?$$
Given a scheme $S$ of finite type over $k$, to each
stable curve $f:C \ra S$ we assign a pseudostable curve $\cT(f):\cT(C)\ra S$.
This recipe must be compatible with base extension: Given a morphism
$h:S_2\ra S_1$ and a family $f:C\ra S_1$, we have a natural
isomorphism
$h^*\cT(C)\stackrel{\sim}{\to}\cT(h^*C)$.
It must also be compatible with isomorphisms: Given two
families $f_1:C_1\ra S$ and $f_2:C_2 \ra S$ and an isomorphism
$\psi:C_1\ra C_2$ over $S$, there is an induced isomorphism over $S$
$$\cT(\psi): \cT(C_1) \ra \cT(C_2),$$
satisfying the standard compatibilities under composition
and base extension.
(We are just recalling the definition of a $1$-morphism
between two categories fibered in groupoids over the category of schemes
cf. \cite{LM}, ch. 2.)

Our main tool is a replacement morphism $\xi_C:C\ra \cT(C)$ defined
over $S$, with the
following properties:
\begin{enumerate}
\item[(a)] over the open subset $S_0 \subset S$ mapping
to the complement to $\delta_1$,  $\xi_C$ is an isomorphism;
\item[(b)] for $s\in S$ mapped to $\delta_1$,
$\xi_{C_s}:C_s\ra \cT(C)_s$ replaces each elliptic tail in $C_s$
with a cusp;
\item[(c)] $\xi_C$ is compatible with base extension
and isomorphisms as indicated above.
\end{enumerate}

We shall specify $\cT$ and $\xi$ using a suitable atlas of $\FMg$ in the
smooth or flat topology.  Let $\pi:Z\ra U$ be a stable curve over a
scheme so that the classifying map $\rho_\pi:U\ra \FMg$ is
faithfully flat. For example, we could take $U$ to be $\Ch^s_3$ (or
the corresponding Hilbert scheme) and $\pi:Z\ra U$ the universal
family. Let $R\rightrightarrows U$ be the corresponding presentation
of $\FMg$, where $R$ encodes the isomorphisms among the fibers of
$\pi$. To construct $\cT$ and $\xi$ over $\FMg$, it suffices to
construct $\cT(\pi):\cT(Z)\ra U$ and $\xi_Z$ compatibly with the
isomorphism relation $R$. Indeed, given an arbitrary stable curve
$f:C\ra S$ with classifying map $\rho_f:S\ra \FMg$, we can pullback
to the fiber product
\[
\xymatrix{ S':=S\times_{\FMg}U \ar[r]^-{pr_U} \ar[d]_-{pr_S} & U
\ar[d]^-{\rho_\pi} \\
S \ar[r]_-{\rho_f} & \FMg }
\]
Consider the resulting stable curves
\[
pr_S^*f:pr_S^*C\ra S' \quad pr_U^*g:pr_U^*Z\ra S',
\]
with classifying maps $\rho_f\circ pr_S$ and $\rho_\pi\circ pr_U$.
The commutativity of the fiber product diagram implies that our two
stable curves are isomorphic over $S'$. By basechange, we have
\[
pr_U^*\cT(\pi):pr_U^*\cT(Z) \ra S', \quad pr_U^*\xi_Z:pr_U^*Z\ra
pr_U^*\cT(Z)
\]
and thus the corresponding constructions for $pr_S^*f:pr_S^*C\ra
S'$. Now $pr_S$ is the basechange of $\rho_\pi$ and thus is
faithfully flat. Since $\cT$ is compatible with isomorphisms, the
desired
\[
\cT(f):\cT(C)\ra S, \quad \xi_C:C\ra \cT(C)
\]
exist by descent.

\subsection{Sketch construction of the morphism}
Assume $\pi:Z\ra U$ is a stable
curve over a smooth base with faithfully flat classifying map
$\rho_\pi:U\ra \ocM_g$. Let $W=\rho_{\pi}^*\delta_1$.  Let $\mu_\pi
: Z\ra \ocM_{g,1}$ denote the classifying morphism to the universal
family, $\d_{1,\{1\}}\subset \ocM_{g,1}$ the boundary divisor
corresponding to the elliptic tails, and
$E=\mu_\pi^*\delta_{1,\{1\}}$, which is also a Cartier divisor on
$Z$.
\begin{eg} \label{eg:mult}
It may seem counterintuitive that the divisor of elliptic tails
should always be Cartier.  For example, consider a stable curve
$$\pi:Z\ra \Delta=\spec k[[t]]$$
with smooth generic fiber and special fiber
$$\pi^{-1}(0)=C_1\cup_p C_2$$
where $C_1$ is smooth of genus one, $C_2$ is smooth of genus $g-1$,
and $\{p\}=C_1\cap C_2$ is a single node. If $Z$ has local analytic
equation $xy=t^m$ at $p$ then $C_1$ is {\em not} Cartier at $p$
unless $m=1$.  However, note that $E\equiv mC_1$ is a Cartier
divisor.
\end{eg}
Fix
$$L=\omega_\pi(E)$$
where $\omega_\pi$ is the dualizing sheaf for the morphism $\pi$.
Away from $W$ we have $L=\omega_{\pi}$, which yields
\begin{prop} \label{prop:easyfree}
For each $n>0$, $\pi_*L^n$ is locally free of rank
\[
k_n := \begin{cases}
g & \mbox{if }  n = 1,\\
(2n-1)(g-1) &  \mbox{if } n \geq 2,
\end{cases}
\]
over $U\setminus W$.  Furthermore, the higher-direct image
\[
R^1\pi_*L^n \simeq \begin{cases}
\cO_U &  \mbox{if } n=1,\\
0 & \mbox{if } n\ge 2,
\end{cases}
\]
over $U\setminus W$.  In particular,
$$\pi_*L^n|_u  \simeq \Gamma(Z_u,\omega_{Z_u}^n)$$
for each $u \in U\setminus W$.
\end{prop}

The family $\{\pi_*L^n|_u\}_{u \in U}$ of linear systems over $U$
will induce the $U$-morphism $Z \to \mcl T(Z)$.
For the remainder of this section, we focus
our attention on the behavior of these linear systems near a
point $u_0 \in W$.
In \S\ref{sS:localfree} we shall prove the local freeness of
$\pi_*L^n$.  In \S\ref{S:limit}, we give a careful analysis of limits
of linear series $\pi_*L^n|_u$.  This is important for two reasons:
First, the local freeness proof depends on it.  Second, it shows
that each linear system $\pi_*L^n|_u$
replaces the elliptic tails of $Z_u$ with cusps.

\subsection{An explicit complex computing the cohomology of $L$}
Suppose that $Z_{u_0}:=C$ has $r$ elliptic tails $E_1,\ldots,E_r$, and let $D$
denote the union of the remaining components.  In analyzing whether
$\pi_*L^n$ is locally free, we may liberally apply faithfully flat
base-extensions $U'\ra U$.  (Cohomology commutes with flat
basechange, and being locally-free is a local property in the
faithfully-flat topology.) For simplicity, we will not use new
notation for these base extensions.

After such a base extension, we may assume
there exist sections
$s_1,\ldots,s_r:U \ra Z$
so that $s_i(u_0)\in E_i$ as a smooth point of $Z_{u_0}$. Over this
neighborhood, we consider the exact sequence
$$0 \ra L^n \ra L^n(s_1+\ldots+s_r) \ra L^n(s_1+\ldots+s_r)|_{s_1,\ldots,s_r}
\ra 0.$$ The last term is supported in $s_1(u_0),\ldots,s_r(u_0)$, a
finite set, and thus has vanishing higher cohomology. As for the
second term,
$$
H^1(Z_{u_0},L^n(s_1+\ldots+s_r)|Z_u)=H^1(D,L^n|_D)
\oplus
\left(\oplus_{i=1}^r H^1(E_i,\cO_{E_i}(s_i(u_0)))\right)=0.
$$
The resulting long exact sequence takes the form
\[
0 \to \pi_*L^n \to F^0 \strl{\varphi}{\to} F^1 \to  R^1\pi_*L^n \ra 0
\]
with $F^0 := \pi_*(L^n(s_1+\ldots+s_r))$ and $F^1 :=
\pi_*(L^n(s_1+\ldots+s_r)|_{s_1,\ldots,s_r})$ locally free of rank
$r_0$ and $r_1$, respectively.  The cokernel sheaf $Q :={\rm
Coker}(\pi_*L^n \to F^0)$ is a subsheaf of a locally-free sheaf,
hence is locally-free of rank $r_0-k_n$ away from a subset $Y
\subset U$ of codimension $\ge 2$.  Note that $Y$ is contained in
the locus where $R^1\pi_*L^n$ fails to be locally free, and thus is
a subset of $W = \rho_\pi^*\d_1$.

\subsection{Limit of linear series}\label{S:limit}
Most of key ideas of this section can be found in \cite{EH},
which develops a general theory of limiting linear series.
However, the case of elliptic tails behaves especially well;
this will be crucial for our application.

Let $B$ be spectrum of a DVR over $k$, with special point
$0$ and closed point $b$.
Consider a morphism $\beta :
(B,0) \to (U,u_0)$ with $\beta(0)=u_0 \in W$ and
$\beta(b) \not\in W$.   We have the fiber square
\[
\xymatrix{ Z_B \ar[r]^-{\beta'} \ar[d]_-{\pi_B} & Z \ar[d]^-\pi \\
B \ar[r]_-\beta & U}
\]
and we write $L_B=\beta'^*L$.

The direct image ${\pi_B}_*L_B^n$ is locally free of
rank $k_n$--it is a subsheaf of $\beta^*F_0$.
We can think of
\[
 {\pi_B}_*L_B^n|_0 \subset H^0(C, L_B^n|_C) =
H^0(Z_{u_0}, L^n|_{Z_{u_0}})
\]
as  {\it the limiting linear series}
\[
\displaystyle{\lim_{\substack{u \to u_0\\ u \in
B}} H^0(Z_u, L^n|_{Z_u})}.
\]
We shall prove that the limit
does not depend on the curve:
${\pi_B}_*L_B^n|_0$ is
canonically identified with a fixed subspace of $H^0(Z_{u_0},
L^n|_{Z_{u_0})}$, regardless of the choice of $\beta : (B, 0) \to
(U, u_0)$.

\begin{notn} Let $X$ be a complete curve and $M$ be a line bundle on
$X$.  Given a linear series $V \subset \G(X, M)$ and an effective
Cartier divisor $A$, we let $V(-A)$ denote the intersection $V \cap
\G(X, M-A)$ in $\G(X, M)$.
\end{notn}

Note that since the restriction of $L$ to an
elliptic tail $E_i$ is trivial, any element of $H^0(C, L^n|_C)$
that vanishes on $D$ vanishes entirely on $C$ and the restriction
map gives an inclusion of $H^0(C, L^n|_C)$ into
$$H^0(D, L^n|_D) =
H^0(D, \o_D^n(2n \sum_{j=1}^r p_j)).$$
In particular, we obtain an inclusion
$$V_n(-a\pj)\subset
H^0(D, \o_D^n((2n-a) \sum_{j=1}^r p_j)).$$

\begin{prop}\label{P:V_n}
Retain the notation introduced above.
For $n \ge 1$, $V_n := {\pi_B}_*L_B^n|_0$ is
naturally identified with
\[
\G(D, \o_D^n((2n-2)\sum_{j=1}^r p_j)) + \lan \sm_1^n, \dots,
\sm_r^n \ran
\]
as a subspace of $H^0(C, L^n|_C)$, where $\sm_j$ are sections of
$V_1$ such that $\sm_j(p_j) \ne 0$ and $\sm_j$ vanishes to order two
at $p_i$ for $i \ne j$.  In particular, $V_n$ satisfies the
following properties:

\begin{enumerate}
\item The vanishing sequence of $V_n$ (regarded as linear series on $D$)
at each $p_j$ is
\[
(0, 2, 3, \dots).
\]
\item
For $n\ge 2$ and $2\le a \le 2n-2$, $V_n$ has $r$ sections
$\t^{(a)}_{n,j}$, $j = 1, \dots, r$, whose images form a basis for
the quotient space
\[
V_n(-a \pj)/V_n(-(a+1)\pj).
\]
\end{enumerate}
\end{prop}

\begin{coro} \label{coro:tocusps}
For $n\ge 2$, the linear series $V_n$ defines a morphism
$$C \ra  \bP^{k_n-1}$$
with image $\cT(C)$.  The induced $C\ra \cT(C)$ is the
morphism introduced in Proposition~\ref{prop:replacing}.
\end{coro}
\begin{proof}[Proof of the Corollary, assuming Proposition~\ref{P:V_n}]
The linear series $V_n$ contains
$$\G(D, \o_D^n
((2n-2)\sum_{j=1}^r p_j)).$$
For any pointed stable curve $(D,p_1,\ldots,p_r)$,
$(\omega_D(p_1+\ldots+p_r))^2$ is very ample
unless $g(D)=2$ and $r=0$, which cannot occur.
Thus the sections of $V_n$ induce an imbedding away on $D$
from $p_1,\ldots,p_r$.  The vanishing sequences imply that
these points are mapped to distinct cusps.  Finally, sections
in $V_n$ are constant along the $E_i$, but not identically
zero, so each $E_i$ is collapsed to the corresponding cusp.
\end{proof}

\begin{proof}[Proof of Proposition~\ref{P:V_n}]
Our whole analysis is over $B$;
for notational simplicity we drop subscripts,
i.e., $L$ designates $L_B$, $\pi$ designates $\pi_B$, and
$E$ designates $E_B$.

For each $n\ge 1$, we have the exact sequence
$$0 \ra \omega^n_{\pi}((n-1)E) \ra L^n \ra L^n \otimes \SS_E \ra 0$$
and on taking cohomology, we obtain
$$0 \ra \pi_*\omega^n_{\pi}((n-1)E) \ra \pi_* L^n \ra
\pi_*(L^n \otimes \SS_E).$$
The cohomology of the first term is locally free
and commutes with base extension.
Indeed, $\omega^n_{\pi}((n-1)E)$ is the dualizing sheaf for $n=1$
and has vanishing higher cohomology for $n>1$.  Thus we deduce
$$\pi_*\omega^n_{\pi}((n-1)E)|_0\simeq
\Gamma(C,\omega^n_{\pi}((n-1)E)|_C),$$
where the latter term is the kernel of
$$\gamma:\Gamma(D,\omega^n_{\pi}((n-1)E)|_D)\oplus
\left(\oplus_{i=1}^r \Gamma(E_i,\omega^n_{\pi}((n-1)E)|_{E_i})\right) \ra
\oplus_{i=1}^r \omega^n_{\pi}((n-1)E)|_{p_i}.$$
This is the gluing condition for the partial normalization of $C$
at the nodes $\{p_1,\dots,p_r\}$.  More concretely, we have
$$\gamma:\Gamma(D,\omega^n_D((2n-1)(p_1+\ldots+p_r)))\oplus
\left(\oplus_{i=1}^r \Gamma(E_i,\omega_{E_i}(p_i))\right) \ra k^r.$$
Since $p_i$ is a basepoint of $\omega_{E_i}(p_i)$, we find
$$\Gamma(C,\omega^n_{\pi}((n-1)E)|_C)\simeq
\Gamma(D,\omega^n_D((2n-2)(p_1+\ldots+p_r)))\oplus
\left(\oplus_{i=1}^r \Gamma(E_i,\omega_{E_i})\right).$$

Our next task is to analyze the image of
$$\eta:\pi_*\omega^n_{\pi}((n-1)E)|_0 \ra \pi_*L^n|_0.$$
Since this is induced by multiplication by $E$, the sections coming
from $\eta$ necessarily vanish along $E$.  Thus $\eta$ naturally
factors through
$$\Gamma(D,L^n|_D\otimes \SS_D(-p_1-\dots-p_r))\simeq
\Gamma(D,\omega_D((2n-1)(p_1+\dots+p_r))).$$
In light of the analysis above, we find
$$\mathrm{image}(\eta)=\Gamma(D,\omega_D((2n-2)(p_1+\ldots+p_r)))\subset V_n.$$

We produce sections $\sm_i\in V_1,i=1,\dots,r$,
which do not vanish at $p_i$ and vanish to order two at $p_j$
for $j \ne i$.
Write $E=\sum_i m_iE_i$ so that each $m_iE_i$ is Cartier
and meets $D$ with multiplicity one (cf. Example~\ref{eg:mult}).
Consider the exact sequence
\[
0 \ra \o_\pi(m_i E_i) \ra L \ra L\otimes \SS_{E\setminus E_i} \ra
0
\]
and the induced map $\iota : \pi_*(\o_\pi(m_iE_i))|_0 \ra V_1$.
It sufficess to find a section
$s_i \in \pi_*(\o_\pi(E_i))|_0$
that does not vanish at $p_i$. Let $E_i^c$ denote the closure of the
complement of $E_i$.
\begin{enumerate}
\item Since $\o_\pi(m_iE_i)|_{E_i} = \SS_{E_i}$, we have
$\pi_*(\o_\pi(m_iE_i))|_0 \subset \G(E_i^c, \o_{E_i^c}(2p_i))$.

\item Since $\G(E_i^c, \o_{E_i^c}) = \G(E_i^c, \o_{E_i^c}(p_i))$, if $p_i$
is a base point of $\pi_*(\o_\pi(m_iE_i))|_0$, then
$\pi_*(\o_\pi(m_iE_i))|_0 \subset \G(E_i^c, \o_{E_i^c})$, which is
absurd from dimensional considerations.
\end{enumerate}
Hence there is a section $s_i \in \pi_*(\o_\pi(m_iE_i))|_0$ that does
not vanish at $p_i$.  The map $\iota$ fits in the following
commutative diagram
\[
\xymatrix{  \pi_*(\o_\pi(m_iE_i))|_0 \ar[r]^-{\iota} \ar@{^(->}[d] &
V_1
\ar@{^(->}[d]\\
\G(E_i^c, \o_{E_i^c}(2p_i)) \ar[r]_-{rest} & \G(D,
\o_{D}(2\sum_{j\ne i}^r p_j)) }
\]
from which it is clear that $\iota(s_i) \ne 0$. Moreover, in view of
the restriction map $rest$, $\sm_i := \iota(s_i)$ vanishes to order
two at $p_j$, $j \ne i$.  In particular, the images of $\sm_j$'s in
the quotient space $\G(D, \o_D(2 \sum_{j=1}^r p_j))/\G(D,
\o_{D})$ are linearly independent. Hence $V_1 \supset \G(D, \o_{D})
+ \lan \sm_i, \dots, \sm_r \ran$.  But both sides are of dimension
$g$ and the equality follows.

Consider the multiplication map $\Sym^n V_1 \to
V_n$ induced from
\[
\Sym^n \pi_*L \to \pi_*L^n.
\]
Thus $V_n$ contains
$\sm_1^n, \dots, \sm_r^n$.
Because of the vanishing properties of the $\sm_j^n$'s, they are
linearly independent in the quotient space $V_n/\G(D, \o_{D}^{\ten
n}((2n-2)(p_1+\dots+p_r)))$. By comparing the dimensions,  we conclude that
$$V_n =
\G(D, \o_{D}^{\ten n}((2n-2)(p_1+\dots+p_r)) + \lan \sm_1^n, \dots, \sm_r^n\ran.$$

Finally, we produce the desired sections $\t^{(a)}_{n,j}$.
Observe that for $a\ge 2$
$$
V_n(-a \pj)=\Gamma(D,\omega^n_D((2n-a)\pj))
$$
so Riemann-Roch and the vanishing of higher cohomology gives
$$\dim V_n(-a\pj)=2n(g-r-1)+(2n-a)r$$
for $2\le a\le 2n-1$.
Thus the quotient
$$V_n(-a\pj)/V_n(-(a+1)\pj)$$
has dimension $r$ provided $2\le a\le 2n-2$.  A suitable basis
lifts to the $\t^{(a)}_{n,j}$.
\end{proof}

\subsection{Local freeness of $\pi_*L^n$}\label{sS:localfree}

\begin{prop}\label{P:makeT1}
For each integer $n\ge 1$, $\pi_*L^n$ is locally free of rank $k_n$.

\end{prop}
\begin{proof}

Given a locally-free sheaf $F$ on $U$ of rank $r$, let $\Gr(m,F)$
denote the scheme representing the functor
\[
\FGr(m, F) : \mbox{ $U$-schemes} \to {\rm Sets}
\]
that associates to a $U$-scheme $\beta : B \to U$ the set of
rank-$m$ subbundles of $\b^*F$, or equivalently, rank-$(r-m)$
locally-free quotients of $\b^*F$.  Our map $\varphi : F^0 \to F^1$
induces a morphism
\[
\tau: U\setminus Y \to \Gr(r_0-k_n, (F^0)^*) \ex \Gr(r_0-k_n,F^1)
\]
and a rational map over $Y$.
Via the Pl\"ucker and Segre embeddings
\[
\begin{array}{l}
 \Gr(r_0-k_n, (F^0)^*) \ex \Gr(r_0-k_n, F^1) \inj
\P(\bigwedge^{r_0-k_n}(F^0)^*) \ex \P(\bigwedge^{r_0-k_n}F^1)  \\
\inj \P(\bigwedge^{r_0-k_n}(F^0)^*\ten \bigwedge^{r_0-k_n}F^1) =
\P(\mathcal{H}om(\bigwedge^{r_0-k_n}F^0, \bigwedge^{r_0-k_n}F^1))
\end{array}
\]
we consider $\tau$ as a rational map  $U  \dra \P(\mcl
Hom(\bigwedge^{r_0-k_n}F^0, \bigwedge^{r_0-k_n}F^1))$ given by the
section $\wedge^{r_0-k_n}\varphi$. The indeterminacy locus of $\tau$
is precisely the zero locus of $\wedge^{r_0-k_n}\varphi$ which is
defined by the ideal generated by $(r_0-k_n+1) \ex (r_0-k_n+1)$
minors of $\varphi$. This is the first Fitting ideal $\mcl I$ of
$R^1\pi_*L^n$. Blowing up along $\mcl I$ to resolve indeterminacy,
we get
\[
\xymatrix{ \til{U}:= Bl_\mcl I(U) \ar[dr] \ar[d]_-{\sm} \\
U \ar@{-->}[r] & \Gr(r_0-k_n, (F^0)^*)\ex\Gr(r_0-k_n, F^1). }
\]
Let $\wt{\tau}$ denote the lift of $\tau$ to $\til{U}$.
We claim that $\sm$ is an isomorphism so that $\wt{\tau}$ is defined
on $U$ and extends $\tau$. This implies that there exists a {\it
subbundle} $K$ of $F^0$ over $U$ such that $K|_{U\setminus Y} \simeq
\pi_*L^n|_{U\setminus Y}$.  For a general discussion of flattening
and blowing up, we refer the reader to \cite{Raynaud} chapter 4.

Since $\sm$ is proper birational and $U$ is normal, it suffices to
show that $\sm$ is quasi-finite. Suppose that it is not quasi-finite
and there exist two distinct closed points $x_1, x_2 \in \til{U}$
such that $\tau$ is not regular at $u:=\sm(x_1) = \sm(x_2)$. Choose
curves $\til{\beta_i} : B = \spec R \to \til{U}$ with
$\til{\beta}_i(0) = x_i$ and $\til{\beta}_i(r) \not\in \d_1$, $i =
1, 2$. Here $R = k[[t]]$,  $0$ denotes the closed point and $r$, the
generic point. We can assume that $\til\tau\circ \til{\beta}_1(0)
\ne \til\tau\circ\til{\beta}_2(0)$ since otherwise $\til\tau$
descends to a (regular) morphism on $U$ at $u$.  Let $\beta_i =
\sm\circ\til{\beta}_i$. Due to the functoriality of $\Gr(k_n, F^0)$,
$\til\tau\circ\til{\beta_i}$ corresponds to a locally-free quotient
$0 \to K_i \to \beta_i^*F^0 \to Q_i \to 0$ such that $\rk(K_i) =
k_n$. Consider the fibre square:
\[
\xymatrix{ \ar@{}[dr] |{\square} Z_B \ar[r]^-{\mu_i} \ar[d]_-{\pi_i}
& Z
\ar[d]^-{\pi} \\
B \ar[r]_-{\beta_i} & U}
\]
We have the exact sequence
\[
0 \lra \pi_{i*}\mu_i^*L^n \lra \beta_i^*F^0 \strl{\beta_i^*\varphi}{\lra}
\beta_i^*F^1 \lra R^1(\pi_i)_*\mu_i^*L^n \lra 0.
\]
Note that the quotient sheaf $\beta_i^*F^0/\pi_{i*}\mu_i^*L^n$ is
locally free, and $\pi_{i*}\mu_i^*L^n$ and $K_i$ are subbundles of
$\beta_i^*F^0$ that agree on $B\setminus{\{0\}}$. It follows that
they are isomorphic on $B$.

We proved in the previous section (Proposition~\ref{P:V_n}) that
$[(\pi_1)_*\mu_1^*L^n]|_0$ and  $[(\pi_2)_*\mu_2^*L^n]|_0$ are
identified in $ H^0(Z|_{\beta_i(0)}, L^n|Z|_{\beta_i(0)})$. This
implies that $\til\tau\circ\til{\beta}_1(0) =
\til\tau\circ\til{\beta}_2(0)$ and that  there exists a {\it
subbundle} $K$ of $F^0$ over $U$ such that $K|_{U\setminus Y} \simeq
\pi_*L^n|_{U\setminus Y}$. Let $j$ denote the natural inclusion
$U\setminus Y \inj U$. Since $\pi_*L^n$ is the kernel of a
homomorphism of locally-free sheaves, it is reflexive and
$j_*\left(\pi_*L^n|_{U\setminus Y}\right) \simeq \pi_*L^n$ since
$\mathrm{codim}(Y)\geq 2$. Therefore, we have
\[
K \simeq j_*\left(K|_{U\setminus Y}\right) =
j_*\left(\pi_*L^n|_{U\setminus Y}\right) \simeq \pi_*L^n
\]
and hence $\pi_*L^n$ is a vector subbundle of $F^0$.
\end{proof}

\subsection{Analysis of the morphisms}
\begin{prop} \label{P:makeT2}
For $n\ge 2$, the sections of $L^n$ relative to $\pi$ induce a
morphism $Z\ra \P(\pi_*L^n)$ over $U$, which factors
$$Z\stackrel{\xi_Z}{\ra} \cT(Z) \hookrightarrow \P(\pi_*L^n)$$
where the second arrow is a closed embedding.
\end{prop}
\begin{proof} By Proposition~\ref{P:makeT1} $\pi_*L^n$ is
locally free.  For $u\in U\setminus W$, Proposition~\ref{prop:easyfree}
guarantees $\pi_*L^n|_u=
\Gamma(\omega_{Z_u}^{\otimes n})$, and $\omega_{Z_u}^{\otimes n}$
is very ample on stable curves of genus $g>2$.  For $u\in W$,
Corollary~\ref{coro:tocusps} implies $\pi_*L^n|_u=V_n$
induces the morphism replacing elliptic tails with cusps.

We therefore have a morphism
$Z\ra \P(\pi_*L^n)$ whose Stein factorization is $\cT(Z)$.
Due to the functoriality of
the relative dualizing sheaf and the line bundle $\d_{1,1}$ on the
moduli stack $\FMgone$, this construction is compatible with
isomorphism relation and commutes with base extension.
\end{proof}

\section{$\Mps$ as a log canonical model for $\Mg$}

In this section, we identify $\Mps$ and $\Mg(9/11)$,
thus completing our proof of Theorem~\ref{T:main1}.
We also prove Theorem~\ref{T:main2}:  A log canonical
model for the pair $(\Mg, K_{\FMg} + \a\d)$ exists as a projective
variety for $7/10 < \a \le 9/11$.

We want to relate the
contraction of the extremal ray $R$ parametrizing elliptic tails
$$\Upsilon: \Mg \ra \Mg(9/11)$$
to the natural birational morphism
$$T:\Mg \ra \Mps.$$
Since $\Mps$ is a GIT quotient, it is a normal projective variety.
Our pointwise description of the natural transformation
$\cT:\FMg \ra \FMps$
in Proposition~\ref{prop:replacing} implies that
\begin{itemize}
\item{$T$ is an isomorphism over $\Mg \setminus \D_1$;}
\item{$T$ takes $R(C_2,p)$ to the cuspidal curve
$\cT(C_1 \cup_p C_2)$, i.e., the curve normalized by
$C_2$ with conductor in $p$;}
\item{$\cT(C)=\cT(C')$ if and only if $C$ and $C'$
have the same number of elliptic tails and
\[
(D,p_1,\ldots,p_r)\simeq
(D',p_1',\ldots,p_r'),
\]
in the notation of Proposition~\ref{prop:replacing}.}
\end{itemize}
In particular, the only curves contracted by $T$ are in the class
$R$, i.e., $T$ is also an extremal contraction of $R$.
Uniqueness of extremal contractions \cite{KMCC} 1.26
implies $T=\Upsilon$ and $\Mps\simeq \Mg(9/11)$,
and the proof of Theorem~\ref{T:main1} is complete.

We turn to the proof of Theorem~\ref{T:main2}.
Let $\dps$ denote the boundary divisor of $\FMps$;
it is the image of $\delta$ under $\cT$.
The proper transform of $\dps$ in $\FMg$ is
$\d-\d_1$.
\begin{lemma}[Log discrepancy formula]
\[
K_{\FMg} + \a \d =
T^*(K_{\FMps} + \a \dps) +
  (9 - 11\a) \d_1
\]
\end{lemma}
\begin{proof}
We determine the value of the discrepancy $c$ in
\[
K_{\FMg} + \a \d =
T^*(K_{\FMps} + \a \dps) +
  c \d_1
\]
by intersecting both sides with the contracted extremal ray $R$.
We have
\[
\begin{array}{ccl}
(K_{\FMg}+\a \d).R&=&
(13\l - 2 \d + \a \d ).R = 13 \l.R - (2 - \a)(\d_0.R +
\d_1.R) \\
&=& 13/12 - (2 - \a)(1 - 1/12)=(11\a-9)/12
\end{array}
\]
and
$c\d_1.R=-c/12$;
the term $T^*(K_{\FMps}+\a \dps)$ obviously yields zero.
\end{proof}

Thus $(\FMps,\alpha \dps)$ has log canonical singularities
for $\a\le 9/11$ provided
$(\FMg,\alpha \d - (9-11\a)\delta_1)$ has log canonical singularities.
The argument
sketched in the Introduction still applies away from $\Delta_1$.
However, reducing the coefficient of $\d_1$ can only increase the
discrepancies of divisors lying over $\d_1$
\cite{FA} 2.17.3;  this does not affect
whether the singularities are log canonical.

To analyze whether $K_{\FMps}+\alpha \dps$ is ample
for $7/10<\a < 9/11$, we describe the ample
cone of $\Mps$ in terms of the ample cone of $\Mg$.  We have
$$T^*\mathrm{NS}(\Mps)=R^{\perp} \subset \mathrm{NS}(\Mg),$$
the hyperplane spanned by a facet of the nef cone
of $\Mg$.  The interior of this facet corresponds to the
ample divisors of $\Mps$.  It suffices then to check:
\begin{prop}  \label{prop:ampleps}
For $7/10 <\a \le 9/11$, the divisor
\[
K_{\FMg} + \a \d -(9-11\a)\d_1
\]
lies in the interior of the facet $R^{\perp}$ of the nef
cone of $\Mg$.
\end{prop}

The following conjectural description of the ample cone of $\Mgn$ would end the
discussion in one stroke:
\begin{conj}[Fulton's conjecture \cite{GKM}] Every one-dimensional
facet of the closed cone of effective curves
$\bar{NE}_1(\bar{M}_{g,n})$ is generated by a one-dimensional boundary
stratum.  Equivalently, any effective curve in $\bar{M}_{g,n}$ is an
effective combination of one-dimensional strata.
\end{conj}
We recapitulate the description of
the one-dimensional strata in \cite{Faber} and
\cite{GKM}:  Let $X_0$ be a 4-pointed
stable curve of genus zero with one point moving and the other three
fixed.
\begin{enumerate}
\item[(a)]{
a family of elliptic tails;}
\item[(b)]{attach a fixed 4-pointed curve of genus $g-3$ to $X_0$;}
\item[(c)]{attach a fixed pointed curve of genus $i$ and a
3-pointed curve of genus $g-2-i$ to $X_0$;}
\item[(d)]{attach fixed two 2-pointed curves of genus $i$ and
$g-2-i$, respectively, to $X_0$;}
\item[(e)]{
attach two 1-pointed curves of genus $i$ and $j$
respectively, and a 2-pointed curves of genus $g-1-i-j$ to $X_0$,
with $i,j \ge 1$ and $i+j = g-1$;}
\item[(f)]{
attach four pointed curves of genus $i$,$j$,$k$, and
$l$, respectively, to $X_0$, with $i+j+k+l = g$ and $i,j,k \ge 1$.}
\end{enumerate}
The intersection of a divisor $a \l - \sum_{i=0}^{\lfloor g/2
\rfloor} b_i \d_i$ with the one-dimensional strata is as follows:
\begin{center}
\begin{tabular}{|c|c|}
  \hline
 (a) & $a-12b_0+b_1$ \\ \hline
(b) & $b_0$ \\
\hline (c) & $b_i$ \\ \hline
 (d) & $2 b_0 - b_i$\\ \hline
 (e) & $b_i + b_j - b_{i+j}$ \\ \hline
 (f) & $b_i + b_j + b_k + b_l - b_{i+j} - b_{i+k} - b_{i+l}$ \\
  \hline
 \end{tabular}
  \end{center}
Here we use the convention that $b_i=b_{g-i}$ for $i>g/2$.

Fortunately, we do not need the full strength of the conjecture.
The following special case is sufficient:
\begin{prop}\label{P:GKM}(6.1,  \cite{GKM})
Let $D = a \l - \sum_{0\le i \le g/2}
b_i\d_i$ be a divisor on $\Mg$ such that either $b_i = 0$ or $b_i
\ge b_0$, $ 1 \le i \le \lfloor g/2 \rfloor$. If $D$ has
non-negative intersection with all the one-dimensional strata, then
$D$ is nef.
\end{prop}

In our case, we have
\[
a = 13, \quad  b_i =
\begin{cases} 2 - \a, \quad i \ne 1 \\
11 - 12 \a, \quad i = 1
\end{cases}
\]
and our divisor satisfies the
hypothesis of Proposition~\ref{P:GKM} as long as $\a \le 9/11$.
Moreover, it
intersects positively with all one-dimensional strata except (a)
as long as the following inequalities hold:
\begin{itemize}
\item{from (b), (c) and (f), $2 - \a > 0$ which is already
satisfied;}
\item{ from (d) with $i=1$ and (e) with $i + j = g-1$,
 $2(2-\a) - (11-12\a)
> 0$ i.e. $\a > 7/10$.}
\end{itemize}
Proposition~\ref{prop:ampleps} therefore follows from Proposition~\ref{P:GKM}.
This completes the proof of Theorem~\ref{T:main2}.

\begin{rmk} \label{rmk:codimtwo}
$T^*(K_{\FMps} + 7/10\ \dps)$ contracts $\D_1$ and the following
loci:
\begin{itemize}
\item{ the stratum (d) with $i = 1$ induces the contraction of
\[
T_0 = \{ C_1 \cup_{p,q} C_2 \, | \, g(C_1) = 1, g(C_2) = g-2\},
\]
the locus of two-pointed elliptic tails or elliptic bridges;}
\item{ the stratum (e) with given $i+j=g-1$ induces the
contraction of
\[T_i = \{ C_1\cup_p C_2 \cup_q C_3 \, | \, g(C_1) =
i, g(C_2) = 1, g(C_3) = g-1-i \} .
\]
}
\end{itemize}
These have codimension two in the moduli space;  at
$\a=7/10$ we obtain a small contraction.
\end{rmk}
\begin{rmk} For special values of $\a \in (7/10,9/11]$,
we can prove the ampleness of $K_{\FMps} + \a\dps$ using invariant
theory.  Recall from \S\ref{S:schubert} that
$$\Mps \simeq \Ch_n \mod \SL_{(2n-1)(g-1)}, \quad n=3,4,$$
where the Chow variety has its natural linearization.
By Theorem 5.15 of \cite{M}, the corresponding line bundle
on $\Mg$ equals $n(g-1)(n(12\l-\d) - 4\l)$,
which is proportional to
$11\l - \d = K_{\FMg} + 9/11 \d$ for $n = 4$, and to $32\l - 3\d
\sim K_{\FMg} + 25/32 \, \d$ for $n = 3$.  These are not
ample of $\Mg$; $\Delta_1$ is a fixed component.
However, the corresponding divisors $K_{\FMps} + \a\dps$ on $\Mps$
are ample for $\a = 9/11, 25/32$.
\end{rmk}

\appendix
\section{Log canonical models of stacks}
\label{sect:LCMS}
Here are conventions for this section:  A `scheme' is a
separated scheme of finite type over a base field.
A `stack' is a Deligne-Mumford stack, separated and of finite type
over a base field;  a `morphism' is a 1-morphism over the base field.
For a general discussion of local properties of representable
morphisms, we refer the reader to \cite{LM} 3.10 and \cite{EGAIVb} 2.6 and 2.7;
we use these without specific attribution.  Roughly,
every property that behaves well under \'etale basechange makes
sense for Deligne-Mumford stacks.
\begin{defn}
A {\em birational proper morphism} of stacks is a proper
representable morphism $f:\cX_1\ra \cX_2$ such that there exist
dense open substacks $\cU_1\subset \cX_1$ and $\cU_2\subset \cX_2$
with $\cU_2=f^{-1}(\cU_1)$ and $f:\mcl U_1 \ra \cU_2$ an
isomorphism. We say that $\cX_1$ and $\cX_2$ are {\em properly
birational} if there exists a stack $\cY$ and birational proper
morphisms $g_1:\cY\ra \cX_1$ and $g_2:\cY \ra \cX_2$.
\end{defn}
It is straightforward to check that this is an equivalence relation.
There is a distinguished open substack $\cU\subset \cX_1$ which
is the largest open set over which $f$ is an isomorphism;
its complement is called the {\em exceptional locus} $\Ex(f)$.  For
any closed substack $\cD\subset \cX_2$ with $f(\cU)\cap \cD$ dense in
$\cD$, the {\em birational transform}
$f_*^{-1}\cD$ is defined as the closure of $f^{-1}(\cD\cap f(\cU))$
in $\cX_1$.
\begin{prop}
\label{prop:resolve} Let $\cX$ be a stack over a field of
characteristic zero.  Then there exists a smooth stack
$\widetilde{\cX}$ and a birational proper morphism
$f:\widetilde{\cX} \ra \cX$; this is called a {\em resolution} of
$\cX$. Furthermore, $\Ex(f)$ can be taken to be a normal crossings
divisor in $\widetilde{X}$. If $\cZ\hookrightarrow\cX$ is a closed
substack with ideal sheaf $I_{\cZ} \subset \cO_{\cX}$ then $f$ can be
chosen so that $f^*I_{\cZ}$ is an invertible sheaf
$\cO_{\tilde{\cX}}(-\widetilde{\cD})$ and $\Ex(f)\cup
\widetilde{\cD}$ is simple normal crossings. This is called a {\em
log resolution of $\cZ$ in $\cX$.}
\end{prop}
\begin{proof}
This statement should be compared to
Theorem 0.2 of \cite{KMCC}.  It is an application of the functorial
procedure for resolving singularities \cite{BM}:  Such procedures
commute with \'etale base extension \cite{Hau}, pp. 329.
In particular, given a \'etale presentation of $\cX$
$$R \rightrightarrows U$$
canonical resolutions $\widetilde{R}\ra R$ and $\widetilde{U}\ra U$
form a presentation
$$\widetilde{R} \rightrightarrows \widetilde{U}$$
for a smooth stack $\widetilde{\cX}$.  The induced representable
morphism $\widetilde{\cX} \ra \cX$ is birational and proper. If
$E\subset U$ denotes the closed subscheme corresponding to $\cZ$
then the canonical embedded resolution of $E$ in $U$ naturally
induces a resolution of $\cZ$ in $\cX$.
\end{proof}
\begin{rem}
Dan Abramovich has pointed out that
the 1992 Harvard thesis of Andrew Joel Schwartz,
{\em Functorial Smoothing Of Morphisms In Equal Characteristic 0 },
should also suffice for our purposes.
\end{rem}
Let $X$ be a normal connected scheme, separated and of finite type
over a field, as usual. Let $D=\sum a_j D_j$ be a $\bQ$-divisor on
$X$, with the $D_j\subset X$  distinct codimension-one reduced
closed subschemes;  we assume each $a_j \in [0,1]$.
\begin{prop} \label{prop:local1}
The following properties of $(X,D)$ are local in the \'etale
topology on $X$:
\begin{enumerate}
\item{$X$ is normal;}
\item{$D_j$ is a codimension-one reduced closed subscheme;}
\item{$m(K_X+D)$ is Cartier for some $m>0$.}
\end{enumerate}
A pair $(X,D)$ satisfying the first three conditions is said to be
{\em admissible}.
Given a coherent $\cO_X$-module $F$, the following properties
are local in \'etale topology:
\begin{enumerate}
\item{$F$ satisfies Serre's condition $S_k$ for some fixed $k>0$;}
\item{$F$ is locally free.}
\end{enumerate}
\end{prop}
\begin{proof}
These are standard properties
of descent:  See \cite{EGAIVb} \S 6.1, 6.4, 6.5, and 6.8 for the behavior
of dimension, property $S_k$, property $R_k$, and other standard singularity
conditions under faithfully-flat basechange.  See \cite{EGAIVb} 2.5.1 for
the stability of local-freeness under faithfully-flat base extension.
\end{proof}
\begin{prop}\label{prop:reflex}
Let $X$ be a normal integral scheme
with field of fractions $K(X)$,
and $D$ an integral Weil divisor on $X$.  Let $\cO_X(D)$ denote
the $\cO_X$-module associated to
$$\{f\in K(X): \mathrm{div}(f)+D \ge 0\}.$$
This is a coherent $\cO_X$-module of rank one, satisfying
Serre's condition $S_2$.  The formation of $\cO_X(D)$
commutes with \'etale base-extension.
\end{prop}
\begin{proof}
Let $j:U\subset X$ be the locus
where $X$ is smooth, which is compatible with \'etale
basechange.  Then $\cO_X(D)=j_*\cO_U(D\cap U)$, which
is coherent, $S_2$, and compatible with base extension.
\end{proof}
Let $\cX$ denote a normal connected stack.
Let $\cD=\sum a_j \cD_j$ be a $\bQ$-divisor on $\cX$,
where the $\cD_j\subset \cX$ are distinct
codimension-one reduced closed substacks;
we assume each
$a_j \in [0,1]$.
\begin{defn}
Suppose $(\cX,\cD)$ is proper and admissible in the sense of
Proposition~\ref{prop:local1}.
Its {\em log canonical ring} is the graded ring
$$R(\cX,\cD):=\oplus_{m\ge 0}
\Gamma(\cX,\cO_{\cX}(mK_{\cX}+\lfloor m \cD \rfloor)),$$
where
$$\lfloor m \cD \rfloor=\sum_j \lfloor m a_j \rfloor \cD_j$$
and the summands are defined via Proposition~\ref{prop:reflex}.
\end{defn}
We refer the reader to \S 2.3 of \cite{KMCC} for definitions of
the various singularities arising in birational geometry.
\begin{prop}[Proposition 5.20 of \cite{KMCC}] \label{prop:local2}
Let $(X,D)$ be an admissible pair, as in Proposition~\ref{prop:local1},
defined over a field of characteristic zero.
The following singularity conditions are local in the
\'etale topology:  terminal, canonical, Kawamata log terminal, or log canonical.
\end{prop}
Hence the following notions are well-defined:
\begin{defn}
$(\cX,\cD)$ is {\em terminal}, (resp.
{\em canonical, Kawamata log terminal,
or log canonical})
if it admits a \'etale presentation with the same property.
It is {\em strictly log canonical} if it is log canonical
and $\cX\setminus \cup_i \cD_i$ is canonical.
\end{defn}
The reader will recall that these notions are defined
in terms of discrepancies \cite{KMCC} \S 2.3, measuring how
the canonical divisor changes under birational morphisms:
\begin{prop}
The admissible pair
$(\cX,\cD)$ is terminal, (resp. canonical, Kawamata log terminal,
or log canonical) if and only if there exists a
log resolution
$f:\widetilde{X} \ra X$
such that
\begin{enumerate}
\item{$\Ex(f)=\cup \cE_i$ is a divisor, $\Ex(f) \cup f^{-1}(\cD)$ is
normal crossings, and $\sum_j f_*^{-1}\cD_j$ is smooth;}
\item{we have a $\bQ$-linear equivalence
\[
K_{\tilde{\cX}}+\sum_j a_jf_*^{-1}\cD_j \equiv f^*(K_{\cX}+\cD)+
\sum_i d(\cE_i;\cX,\cD) \cE_i, \quad d_i:=d(\cE_i;\cX,\cD),
\tag{$\ast$}
\]
with each $a_j<1$ and $d_i>0$ (resp. $a_j\le 1$ and $d_i\ge 0$,
$a_j <1$ and $d_i > -1$, or $a_j\leq 1$ and $d_i\ge -1$.)}
\end{enumerate}
\end{prop}
\begin{proof}
Proposition~\ref{prop:resolve} guarantees there exists a resolution
with the prescribed properties.  By definition, if $(\cX,\cD)$
has the prescribed singularities then the discrepancies $d_i$ and
coefficients $a_j$ satisfying the listed properties.
The converse direction is an application of Corollary 2.32 of \cite{KMCC}.
\end{proof}
\begin{defn}
The {\em discrepancy} $d(\cE_0;\cX,\cD)$ is defined via equation
($\ast$) for any integral codimension-one substack in the
exceptional locus of some birational proper morphism $\cY \ra \cX$.
By convention, we set
$$d_j=d(\cD_j;\cX,\cD)=-a_j$$
and take $d(\cD_0;\cX,\cD)=0$ for any integral codimension-one
substack $\cD_0\subset \cX$ not in the support of $\cD$.
\end{defn}
\begin{defn}\label{defn:propbir}
Two admissible pairs $(\cX,\cD)$ and $(\cX',\cD')$
are {\em properly birational}
if there exists an admissible pair
$(\widetilde{\cX},\widetilde{\cD})$ and
proper birational morphisms $f:\widetilde{\cX} \ra \cX,
f':\widetilde{\cX} \ra \cX'$, so that
$K_{\tilde{\cX}}+\widetilde{\cD}-f^*(K_{\cX}+\cD)$
(resp. $K_{\tilde{\cX}}+\widetilde{\cD}-{f'}^*(K_{\cX'}+\cD')$)
is effective and $f$-exceptional (resp. $f'$-exceptional).
\end{defn}
Our main interest is the case
where $(\cX,\cD)$ and $(\cX',\cD')$ are log canonical
and $(\widetilde{\cX},\widetilde{\cD})$ is a log resolution
with suitably chosen boundary.
\begin{prop}
If the admissible pairs
$(\cX,\cD)$ and $(\cX',\cD')$ are properly birational then
there is a natural isomorphism
$$R(\cX,\cD)\stackrel{\sim}{\ra} R(\cX',\cD')$$
of graded rings.
\end{prop}
\begin{proof}
It suffices to prove the isomorphism
$$R(\cX,\cD)\stackrel{\sim}{\ra} R(\widetilde{\cX},\widetilde{\cD}).$$
We have the discrepancy equation
$$K_{\tilde{\cX}}+\widetilde{\cD}=_{\bQ}f^*(K_{\cX}+\cD)
+\sum_i d_i \cE_i$$ where the $d_i\ge 0$ and the $\cE_i$ are
$f$-exceptional. This was formulated on the level of $\bQ$-Cartier
divisors, but there is a refined interpretation: For each $m\ge 0$,
consider the homomorphisms \smaller
\[
f^*\cO_{\cX}(m K_{\cX}+\lfloor m \cD\rfloor ) \ra
\cO_{\tilde{\cX}}(mK_{\tilde{\cX}}+
    \lfloor m (\widetilde{\cD}-\sum_i d_i \cE_i)\rfloor )
\subset \cO_{\tilde{\cX}}(mK_{\tilde{\cX}}+\lfloor m
\widetilde{\cD}\rfloor). \tag{$\ast\ast$}\] \normalsize
These deserve some explanation: First, there is no harm passing to an
\'etale neighborhood which is a scheme;  for simplicity, we leave
the notation unchanged.  The divisors $\lfloor m f^{-1}_*\cD
\rfloor$ and $\lfloor m (\widetilde{\cX}-\sum_i d_i \cE_i) \rfloor $
agree on all the components which have codimension one in $\cX$, so
we can focus on the $f$-exceptional components.  Every section $s\in
\cO_{\cX}(mK_{\cX}+\lfloor m \cD \rfloor)$ extends to a section with
poles of some order on the exceptional locus of $f$. Replacing $s$
by a suitable power $s^N$, we may assume that $mN(K_{\cX}+\cD)$ is
integral and Cartier, and our discrepancy equation gives the
requisite bound on the pole order. The inclusions $(\ast\ast)$
induce injections
$$\Gamma(m K_{\cX}+\lfloor m \cD\rfloor )\hookrightarrow
\Gamma(mK_{\tilde{\cX}}+\lfloor m \widetilde{\cD}\rfloor ).$$
This is surjective:  Each element of the image is a section of
$$\Gamma(\cX-\cup_i f(\cE_i),m K_{\cX}+\lfloor m \cD\rfloor ),$$
which extends to all of $\cX$ because $\cX$ is normal.
\end{proof}
\begin{prop}\label{prop:redcoarse}
Let $(\cX,\cD)$ be a proper log canonical (resp. Kawamata
log terminal) pair, and let
$\pi:\cX \ra X$ denote the coarse
moduli space of $X$.  Then there exists an effective $\bQ$-divisor
$\check{D}=\sum_{\ell} \check{c}_{\ell} \check{D}_{\ell}, 0\le \check{c}_{\ell} \le 1$
on $X$, with the $\check{D}_{\ell}$ distinct and irreducible
with the following properties:
\begin{enumerate}
\item{
$(X,\check{D})$ is log canonical (resp. Kawamata log terminal);}
\item{for each positive integer $m$ such that
$m(K_X+\check{D})$
is integral and Cartier, there is an equivalence of Cartier divisors
$$m(K_{\cX}+\cD)=\pi^*m(K_X+\check{D});$$}
\item{for every $m\ge 0$ we have
$$\Gamma(X,mK_X+\lfloor m\check{D} \rfloor) \stackrel{\sim}{\ra}
\Gamma(\cX,mK_{\cX}+\lfloor m\cD \rfloor).$$}
\end{enumerate}
Together, these yield a natural isomorphism
$$\pi^*:R(X,\check{D}) \stackrel{\sim}{\ra} R(\cX,\cD).$$
\end{prop}
\begin{proof}
 We first produce the divisor $\check{D}$;  to extract
its coefficients, we localize along each codimension-one integral
closed substack $\cB \subset \cX$.
Let $a(\cB)$ denote the coefficient of $\cD$
at $\cB$, which might be zero if $\cB$ does not appear in $\cD$.
Suppose $\pi$ ramifies to order $e(\cB)$ at the generic point $\eta(\cB)$,
i.e., the order of the automorphism group at the generic
point of $\cB$ is $e(\cB)$ times the order of the automorphism group at
the generic point of $\cX$.  If $B\subset X$ is the reduced image
of $\cB$ then the ramification formula gives
$$\pi^*\omega_{X,\eta(B)}=\omega_{\cX,\eta{\cB}}(-(e-1)\cB) \quad
\pi^*B=e\cB$$
which yields the equation of $\bQ$-Cartier divisors
$$\pi^*(K_X+\sum_{\cB}\frac{e(\cB)-1+a(\cB)}{e(\cB)}B)\equiv K_{\cX}+\cD.$$
We therefore set
$$\check{D}=\sum_{\cB}\frac{e(\cB)-1+a(\cB)}{e(\cB)}B;$$
the pair $(X,\check{D})$ is log canonical (resp. Kawamata log terminal)
by \cite{KMCC} 5.20.

It remains then to check
$$\Gamma(X,mK_X+\sum_B \lfloor m\frac{e(\cB)-1+a(\cB)}{e(\cB)} \rfloor B )
\simeq
\Gamma(\cX,mK_{\cX}+\sum_{\cB} \lfloor m a(\cB) \rfloor \cB).$$
By the ramification formulas, we have
\begin{eqnarray*}
\Gamma(X,mK_X+\sum_B \lfloor m\frac{e(\cB)-1+a(\cB)}{e(\cB)} \rfloor B )
&\hookrightarrow & \\
 \Gamma(\cX,mK_{\cX}+\sum_{\cB}
( e(\cB)\lfloor m\frac{e(\cB)-1+a(\cB)}{e(\cB)}\rfloor-m(e-1)) \cB).
& &
\end{eqnarray*}
This is surjective:
For any $\cO_X$-module $F$ we have $\Gamma(X,F)=\Gamma(\cX,\pi^*\cF)$, as
every section of $\pi^*\cF$ is pulled back from the coarse moduli space.
We must verify that the inclusion
\begin{eqnarray*}
\Gamma(\cX,mK_{\cX}+\sum_{\cB}
( e(\cB)\lfloor m\frac{e(\cB)-1+a(\cB)}{e(\cB)}\rfloor-m(e(\cB)-1)) \cB)
 & \subset   & \\
\Gamma(\cX,mK_{\cX}+\sum_{\cB}\lfloor m a(\cB)\rfloor \cB) & &
\end{eqnarray*}
is an isomorphism.  We analyze sections locally in an \'etale neighborhood
of $\eta(\cB)$, where the stack can be presented as a $\mu_e$-quotient,
with $e=e(\cB)$ (Lemma 2.2.3 of \cite{AbVi});  fix $a=a(\cB)$.    Let $\omega$ be a local
generator of the canonical class and $F$ the local defining equation for
$\cB$, so that $\mu_e$ acts on both $\omega$ and $F$ via multiplication by $\zeta$,
a primitive $e$-th root of unity.  The sheaf
$\cO_{\cX}(mK_{\cX}+\sum_{\cB}\lfloor m a(\cB)\rfloor \cB)$
has local generator
$\omega^{\otimes m}/F^{\lfloor m a \rfloor},$ on which $\mu_e$
acts by $\zeta^{m-\lfloor m a\rfloor}$.  Of course, $\mu_e$ acts
trivially on global sections, so global sections must be multiples of
$$F^{\lfloor m a \rfloor -m - e \lfloor (\lfloor ma \rfloor -m)/e\rfloor}
\omega^{\otimes m}/F^{\lfloor m a \rfloor}=
\omega^{\otimes m}/F^{m + e \lfloor (\lfloor ma \rfloor -m)/e\rfloor}.$$
Thus global sections of
$\cO_{\cX}(mK_{\cX}+\sum_{\cB}\lfloor m a(\cB)\rfloor \cB)$
are forced to have zeros along $\cB$;  precisely, we find
\begin{eqnarray*}
& & \Gamma(\cX,mK_{\cX}+\sum_{\cB}
  (m+e(\cB) \lfloor (\lfloor ma(\cB)\rfloor -m)/e(\cB) \rfloor ))  \cB) \\
& &=\Gamma(\cX,mK_{\cX}+\sum_{\cB}\lfloor m a(\cB)\rfloor \cB).\\
\end{eqnarray*}

The final step is the combinatorial statement
$$e\lfloor m\frac{e-1+a}{e}\rfloor-m(e-1)=
m+e \lfloor ( \lfloor ma \rfloor - m)/e \rfloor,$$
which is equivalent to
$$\lfloor \frac{ma-m}{e}\rfloor=
 \lfloor ( \lfloor ma \rfloor - m)/e \rfloor.$$
Here $e$ and $m$ are positive integers and $a \in [0,1] \cap \bQ$.
If this were false then we could find an integer $n$ with
$$ ( \lfloor ma \rfloor - m)/e  < n \le (ma-m)/e$$
and hence
$$ \lfloor ma \rfloor  < en+m \le ma,$$
violating the definition of the round-down operation.
\end{proof}

The Kawamata basepoint-freeness theorem \cite{KMCC} yields
\begin{coro}[Basepoint-freeness for stacks] \label{coro:bpfs}
Let $(\cX,\cD)$ be a proper Kawamata log terminal pair.
Assume $K_{\cX}+\cD$ is nef and big, i.e., $K_X+\check{D}$
is nef and big, where
$(X,\check{D})$ is the pair constructed in Proposition~\ref{prop:redcoarse}.
Consider the quotient stack and coarse moduli space
$$\cY:=
\left[ (\Spec R(\cX,\cD) \setminus 0)/\bG_m\right] \quad
Y:=
\Proj R(\cX,\cD)=\Proj  R(X,\check{D}),$$
where the action of $\bG_m$ arises from the grading.
Then there is a morphism of stacks
$$\psi:\cX \ra \cY$$
inducing on coarse moduli spaces the contraction from $X$
to the log canonical model of $(X,\check{D})$.
\end{coro}

\bibliographystyle{alpha}
\bibliography{mg911}
\end{document}